\documentclass[12pt]{article}
\usepackage{graphicx}
\usepackage{amsfonts}
\usepackage{subfigure}
\usepackage{amsmath,amssymb}

\newtheorem{theorem}{Theorem}[section]

\newtheorem{corollary}[theorem]{Corollary}
\newtheorem{lemma}[theorem]{Lemma}
\newtheorem{remark}[theorem]{Remark}

\usepackage{hyperref}
\hypersetup{bookmarks=true, bookmarksnumbered=true,
bookmarksopen=true, colorlinks, linkcolor=blue, citecolor=blue,
plainpages=false, pdfstartview=FitH}

\numberwithin{equation}{section}

\begin{document}
\pagestyle{plain}

\title{Orientation-Preservation Conditions on an Iso-parametric FEM in Cavitation
Computation\thanks{The research was supported by the NSFC
projects 11171008 and 11571022.}}

\author{Chunmei Su, \hspace{1mm} Zhiping Li\thanks{Corresponding author,
email: lizp@math.pku.edu.cn} \\  LMAM \& School of Mathematical Sciences, \\
Peking University, Beijing 100871, China}

\date{}
\maketitle
\begin{abstract}
The orientation-preservation condition, i.e., the Jacobian determinant of the
deformation gradient $\det \nabla u$ is required to be positive, is a natural physical
constraint in elasticity as well as in many other  fields. It is well known
that the constraint can often cause serious  difficulties in both theoretical analysis and
numerical computation, especially when the material is subject to large deformations.
In this paper, we derive a set of sufficient and necessary conditions for the quadratic
iso-parametric finite element interpolation functions of cavity solutions to be
orientation preserving on a class of radially symmetric large expansion
accommodating triangulations. The result provides a practical quantitative guide for
meshing in the neighborhood of a cavity and shows that the orientation-preservation can be
achieved with a reasonable number of total degrees of freedom by the quadratic
iso-parametric finite
element method.
\end{abstract}



\section{Introduction}

As early as 1958, Gent and Lindley \cite{Gent and Lindley} carried out physical
experiments and studied the sudden void formation on elastic bodies
under hydrostatic tension. Since then, the phenomenon, which is referred to as
cavitation in literatures, has been intensively studied by numerous researchers.

There are two representative models for the cavity formation. One is the so-called
deficiency model proposed by Gent and Lindley \cite{Gent and Lindley}, in which
the cavities are considered to develop from pre-existing small voids under large
triaxial tensions. The other is the perfect model established by
Ball \cite{Ball82}, in which voids form in an intact body so that the total stored
energy of the material could be minimized. The relations between the two models
are partially established by the work of Sivaloganathan et. al. \cite{Sival2006}
and Henao \cite{Duvan Henao2009}: roughly speaking, given the right positions of the voids,
as the radii of the pre-existing small voids go to zero, the solution of the
deficiency model converges to the solution of the perfect model. Furthermore, the
configurational forces can be used to detect whether a void is formed in the
right position \cite{Lian and Li iso,Sival2002}.

The perfect model typically displays the Lavrentiev phenomenon \cite{Lavrentiev} when
there is a cavitation solution, leading to the failure of the conventional finite
element methods \cite{Bai and Li,Ball and Knowles}.
Though there are existing numerical methods developed to overcome the Lavrentiev
phenomenon (\cite{Bai and Li,Ball and Knowles,Z. Li,Marrero and Betancourt}),
they do not seem to be powerful and
efficient enough to tackle the cavitation problem on their own.

In fact, all of the
numerical studies on cavitation, known to the authors so far, are based on the
deficiency model, in which one considers to minimize the total energy of
the form

\begin{equation}
  \label{functional}
  E(u)=\int_{\Omega_\varrho}W(\nabla u(x))dx,
\end{equation}
in the set of admissible functions
\begin{equation}
  \label{admissible set}
  U=\{u\in W^{1,1}(\Omega_\varrho;\mathbb{R}^n) \ \mbox{is one-to-one
  a.e.}: u|_{\Gamma_0}=u_0, \det \nabla u >0 \ a.e. \},
\end{equation}
where $\Omega_\varrho = \Omega \setminus \bigcup_{i=1}^{K}
B_{\varrho_i}(a_i)\subset\mathbb{R}^{n}\,(n=2,3)$
denotes the region occupied by an elastic body in its reference
configuration, $B_{\varrho_i}(a_i)=\{x\in\mathbb{R}^n:|x-a_i|<\varrho_i\}$ are
the pre-existing defects of radii $\varrho_i$ centered at $a_i$. In
\eqref{functional} $W: M^{n \times n}_+\rightarrow\mathbb{R}^+$ is the stored energy
density function of the material, $M_+^{n \times n}$ denotes the set of $n \times n $
matrices with positive determinant, $\Gamma_0$ is the boundary of $\Omega$. We notice here
that, in elasticity theory, the Jacobian determinant $\det \nabla u$, the local volume
``stretching factor" of a deformation, is naturally required to be positive, which
means that no volume of the material is compressed into a point or even turned ``inside
out". The constraint is of vital importance, for instance it excludes the deformations that
have a reflection component, and it is a necessary condition of the fact that
the matter should not inter-penetrate. On the other hand, the constraint $\det \nabla u>0$,
though less strict than the incompressibility
$\det \nabla u \equiv 1$, also inevitably brings some serious difficulties to
mathematical models (\cite{Ball89, Muller99b}) as well as numerical computations
(\cite{Oscar2015, Lian and Li iso, Xu and Henao}).
\begin{figure}[ht!]

    \centering \subfigure[Failure of affine element.]{

    \includegraphics[width=6.65cm,height=3.8cm]{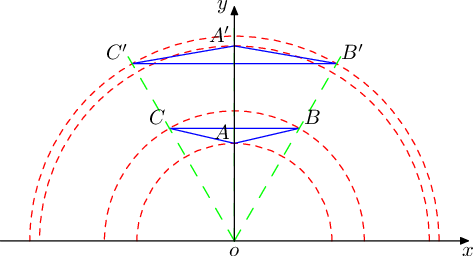}

\label{reverse}} \hspace{3mm} \centering \subfigure[Quadratic iso-parametric
element.]{

  \includegraphics[width=6.65cm,height=3.8cm]{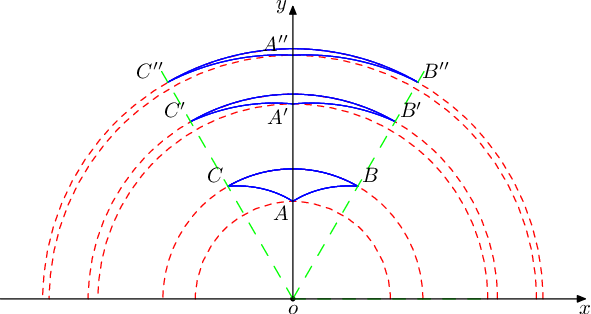}

 \label{noreversion}  }
  \hspace*{-3mm}
  \caption{Quadratic FE is superior in orientation-preservation.} \label{Reversion}

\end{figure}

To have an intuitive view of orientation-preserving behavior of finite element
approximation of cavitation solutions, we compare schematically in Figure~\ref{Reversion}
the affine finite element interpolations and the quadratic
iso-parametric finite element interpolations of a section of a ring before and after a large
radially expansionary deformation. It is clearly seen in Figure~\ref{reverse} that
the affine finite element interpolation fails to preserve the orientation, i.e., the
interpolation triangle $ABC$ before the deformation is anticlockwise but the
interpolation triangle $A'B'C'$ after the deformation is clockwise. It is also easily seen
in Figure~\ref{noreversion} that the quadratic iso-parametric finite element
interpolations can successfully preserve the orientation even for much larger deformations.
This suggests that the conforming affine finite element method is not a good candidate for
the cavitation computation, while quadratic finite element methods might
be. In fact, for the conforming affine finite element method to preserve the orientation,
the amount of degrees of freedom required can be unbearably large \cite{Xu and Henao}.
On the other hand, some numerical methods based on quadratic finite elements
\cite{Lian Dual,Lian and Li iso}, non-conforming affine finite element
\cite{Xu and Henao}, bi-linear and tri-linear finite elements (see \cite{Lopez,Tvergaard}
among many others) have shown considerable numerical success. In particular,
the iso-parametric finite element method developed in \cite{Lian and Li iso} showed great
potential in the computation of multi-voids growth, and the numerical experiments also
revealed that the orientation-preservation conditions are crucial for the method to produce
efficiently accurate finite element cavitation solutions.

The only practical analytical result on the orientation preservation condition for the
cavitation computation known to the authors so far is \cite{SuLiRectan}, where a sufficient condition was given for a dual-parametric bi-quadratic finite element method.

In this paper, we study the orientation-preserving behavior of the quadratic iso-parametric
finite element approximations of cavity solutions by analyzing the sufficient and necessary
conditions for the interpolation functions to preserve the orientation. We will see that,
compared with the dual-parametric bi-quadratic FEM, the derivation of the orientation
preservation conditions for the quadratic iso-parametric FEM is more involved.
Since the cavitation solutions are generally considered to vary mildly except
in a neighborhood of the voids, where the material experiences large expansion
dominant deformations, and where the difficulty of the computation as well as the
analysis lies, we restrict ourselves to a simplified problem with
$\Omega_\varrho = B_1(0) \setminus B_{\varrho}(0)$ in $\mathbb{R}^2$. To bring out the
principal relations in the orientation-preservation conditions and avoid unnecessarily
tedious calculations, we further restrict ourselves to simple expansionary boundary
conditions of the form $u_0=\lambda x$ and the radially symmetric cavitation solutions.
The result shows that the orientation-preservation can be achieved with a reasonable number of
total degrees of freedom. In fact, combined with the corresponding interpolation error
estimates, it would lead to an optimal meshing strategy, which we will show in a
separate article \cite{SuLi}.

The structure of the paper is as follows. In \S~2, we present some properties of the
cavitation solutions for a specific class of energy functionals, as well as
the quadratic iso-parametric
finite element and a radially symmetric large expansion accommodating triangulation method.
\S~3 is devoted to deriving the sufficient and necessary
orientation-preservation conditions on the mesh distribution. We end the paper with some
discussions and conclusion remarks in \S~4.

\section{Preliminaries}

We consider a typical class of stored energy density functions of the form
\begin{equation}
    \label{energy function}
  W(F)=\Phi(v_1,\dots,v_n)=\omega
  \!\left(\sum\limits_{i=1}^n v_i^2\right)^{\frac{p}{2}}+
  g\!\left(\prod_{i=1}^n v_i\right),\quad \forall F \in M^{n\times n}_+,
\end{equation}
where $\omega>0$ is a material constant, $v_1$, $\dots$, $v_n$ are the singular
values of the deformation gradient $F$, and where, to ensure the existence and
regularity of cavity solutions \cite{Muller95}, $n-1<p<n$, and $g:(0,\infty)\rightarrow
[0,\infty)$ is a continuously differentiable strictly convex function satisfying
\begin{equation}\label{asmp g}
g(d)\rightarrow+\infty\mbox{ as } d \rightarrow0,\mbox{ and }
\frac{g(d)}{d}\rightarrow+\infty \mbox{ as } d\rightarrow +\infty.
\end{equation}
For example, $g(d)=\frac{\chi}{2}(d-1)^2+\frac{1}{d}$ was used in \cite{Lian and Li 3}
with the constant $\chi>0$ as the bulk modulus. As mentioned in the introduction, for simplicity, we henceforth assume that $n=2$.

\subsection{Properties of radially symmetric cavitation solutions}

For the simple expansionary boundary condition given by $u_0=\lambda x$, $\lambda>1$,
and the radially symmetric deformations
$u(x)=\frac{r(|x|)}{|x|} x$, the problem
defined on the domain $\Omega_{\varrho} = B_1(0) \setminus B_{\varrho}(0)$
reduces to minimizing the energy of the form
\begin{equation}\label{mini}
I_\varrho(r)=\int_\varrho^1R\Phi\Big(r'(R),\frac{r(R)}{R}\Big) dR
\end{equation}
in the set of admissible functions
$$
A^\lambda_\varrho=\{r \in W^{1,1}(\varrho,1): r(\varrho)> 0,\; r(1)=\lambda,
\; \text{and} \; r'>0 \; a.e. \}.
$$
It is well known (see, {\em e.g.} \cite{Sival1986}) that the problem admits a unique
minimizer $r_\varrho^\lambda \in C^2( (\varrho,1])$, which satisfies the Euler-Lagrange
equation:
\begin{eqnarray}
  \label{EL2d}
  \frac{d}{dR}\left(R\Phi_{,1}\Big(r'(R),\frac{r(R)}{R}\Big)\right)=
  \Phi_{,2}\Big(r'(R),\frac{r(R)}{R}\Big),
  & & R \in (\varrho,1),\\
\label{natural}
\omega p \left(r'(\varrho)^2+\frac{r(\varrho)^2}{\varrho^2}
\right)^{\frac{p}{2}-1} \frac{r'(\varrho)\varrho}{r(\varrho)}+
g'(d(\varrho))=0,\\
\label{Dirichlet}
r(1)=\lambda,
\end{eqnarray}
where $d(\varrho)=\det \nabla u|_{\partial B_\varrho(0)}=
\frac{r(\varrho)}{\varrho}r'(\varrho)$. In particular, for the perfect model
($\varrho=0$), there exists a constant $\lambda_c>1$, such that, for
$\lambda>\lambda_c$, the minimizer satisfies $r_0^\lambda(0)>0$; for
$\lambda\le\lambda_c$, the minimizer is given by $r_0^\lambda(R)=\lambda R$.
Thus, by \cite{Sival1986}, for $\varrho>0$, $r_\varrho^\lambda(R)\ge r_0^\lambda(0)$
when $\lambda >\lambda_c$.

In the case of the perfect model, the
radially symmetric cavity solution $r_0(R)$ is proved to be a bounded strictly convex
function (see \cite{Sival1986}), furthermore, it can be shown that $r_0(R)$
satisfies $m_0 R \le r_0'(R) \le M_0 R$, for all $R\in (0,1]$, where
$0<m_0< M_0$ are constants.
The result is in fact valid also for the deficiency model,
at least when $\varrho$ is sufficiently small.

\begin{lemma}\label{r'r''>0}
Let $r(R)$ be the minimizer of \eqref{mini} over $A^\lambda_\varrho$ with the energy density function given
by \eqref{energy function}. If $\lambda >\lambda_c$, then for sufficiently small
$\varrho\ge 0$, $r(R)$ satisfies
\begin{equation}\label{mR<r'<MR}
0<r''(R) \leq C, \;\;\; m R \le r'(R) \le M R,  \quad \forall R\in [\varrho ,1].
\end{equation}
where $C>0$, $0<m< M$ are constants independent of $\varrho$.
\end{lemma}

\textbf{Proof.}
Firstly, a direct manipulation on \eqref{EL2d} yields
\begin{equation}\label{r''}
r''(R)=\frac{1}{R\Phi_{,11}}\left(\frac{r(R)}{R}-r'(R)\right)
\left(\frac{\Phi_{,2}-\Phi_{,1}}{\frac{r(R)}{R}-r'(R)}+\Phi_{,12}\right).
\end{equation}
It is straightforward to show
that $\Phi_{,11}$ and the term in the second bracket are always positive.
With the same arguments as in \cite{Sival1986}, the term $\frac{r(R)}{R}-r'(R)$
is either identically 0 or never vanishes. If $\frac{r(R)}{R}-r'(R)\equiv 0$,
then $r(R)=\lambda R$, which contradicts the fact that
$r_\varrho^\lambda(R)\rightarrow r_0^\lambda(R)>0$ when $\lambda >\lambda_c$.
On the other hand, it follows from \eqref{natural} that
$g'(d(\varrho))\le 0$. Since $g''(d)>0$, this yields $d(\varrho)=
\frac{r'(\varrho)r(\varrho)}{\varrho} \le d_0$, or equivalently
$r'(\varrho) \le \frac{d_0\varrho}{r(\varrho)}$,
where $d_0$ is the unique root of $g'(x)=0$ (see \eqref{asmp g}). Thus, for
the cavity solution, we have
$\frac{r(\varrho)}{\varrho} - r'(\varrho) \ge \frac{r_c}{\varrho} -
\frac{d_0\varrho}{r_c} >0$, where $r_c=r_0^\lambda(0)$. Hence $\frac{r(R)}{R}-r'(R)>0$
and consequently
$r''(R)>0$, as long as $\varrho < \frac{r_c}{\sqrt{d_0}}$.

Next, we notice that the radial component of the Cauchy stress
$T(R)=\frac{R}{r(R)} \Phi_{,1}($ $r'(R),\frac{r(R)}{R})$ is nowhere
decreasing (\cite{Sival1986}), and $T(\varrho)=0$. Hence $T(R) \ge 0$
for $R \in [\varrho, 1]$,  which can be reformulated as
$$
g'(d(R))\ge - \omega p \left(r'(R)^2+\frac{r(R)^2}{R^2}
\right)^{\frac{p}{2}-1} \frac{r'(R)R}{r(R)}, \quad \forall R \in [\varrho,1).
$$
Since $p\in (1,2)$, $r'(R)R<r(R)$, and $r(R)\geq r_c$, this yields $g'(d(R))\ge
- \omega p r_c^{p-2}$. Thus, by the convexity of $g$ and
$r(R)\leq r(1)=\lambda$, we obtain $r'(R) \geq \frac{d_-}{\lambda} R$, where
$d_- \in (0,d_0)$ is such that $g'(d_-)=- \omega p r_c^{p-2}$.
On the other hand, we notice that $r(R) \geq r_0(R)$ (see \cite{Sival2009}),
and consequently $T_0(R) \geq T(R)$, for all $R \in [\varrho,1)$ and $\varrho>0$,
where $r_0(R)$ is the cavity solution of the perfect model,
$$
T_0(R) = \omega p \left(r_0'(R)^2+\frac{r_0(R)^2}{R^2}
\right)^{\frac{p}{2}-1} \frac{r_0'(R)R}{r_0(R)}+ g'(d_0(R))$$
is the normal surface traction with respect to the perfect model, and
$d_0(R)=\frac{r_0'(R)r_0(R)}{R}$. Thus, we have
$$
g'(d(R))\le g'(d_0(R)) + \omega p \left(r_0'(R)^2+\frac{r_0(R)^2}{R^2}
\right)^{\frac{p}{2}-1} \frac{r_0'(R)R}{ r_0(R)} \leq g'(d_0(R)) +
\frac{\omega p }{r_c^{2-p}},
$$
since $r(R)\geq r_0(R)\geq r_c$ (see \cite{Sival1986}) and $r_0'(R) < \frac{r_0(R)}{R}$.
Denote $d_{0+}= \displaystyle\max_{0 \leq R \leq 1} d_0(R)$ and
$d_+ = (g')^{-1}(g'(d_{0+})
+\omega p r_c^{p-2})$. Then the above inequality yields
$d(R) \leq d_+$, and consequently $r'(R) \leq \frac{d_+}{r_c} R$.
Hence, $m=\frac{d_-}{\lambda}$ and $M=\frac{d_+}{r_c}$ in \eqref{mR<r'<MR}.

The uniform boundedness of
$r''(R)$ can be verified directly by \eqref{r''} using the facts that
$0<d_- \leq d(R) \leq d_+$, $0<r_c \leq r(R) \leq \lambda$ and $mR \leq r'(R) \leq MR$.
\hfill $\square$

\subsection{The quadratic iso-parametric FEM}
\begin{figure}[h!]
 \begin{minipage}[l]{0.5\textwidth}
 \centering \includegraphics[width=2in]{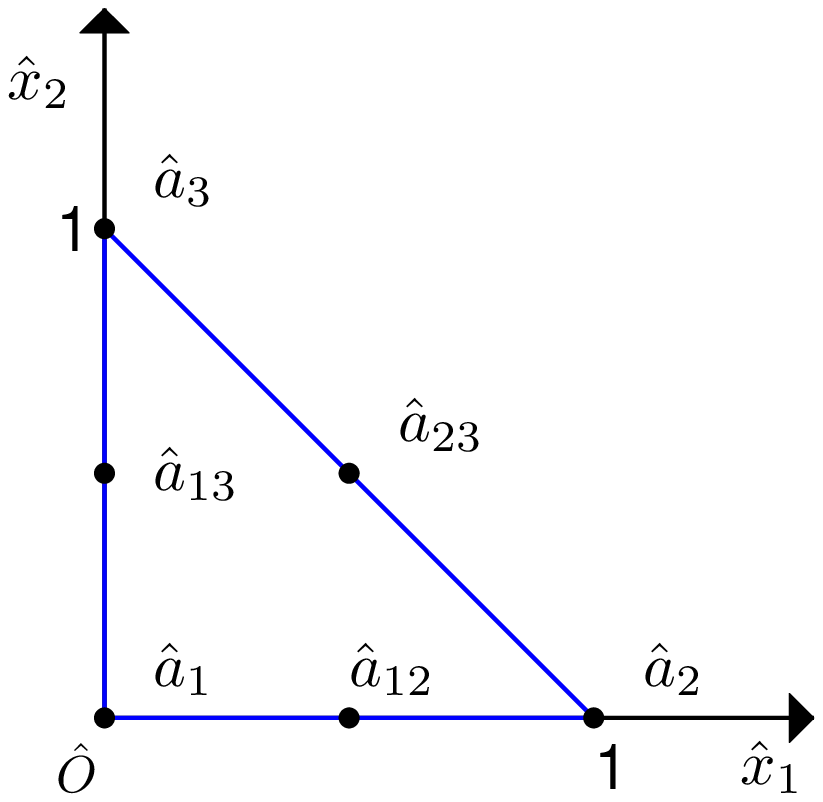}
 \caption{The reference element $\hat{T}$.} \label{dengcan reference element}
\end{minipage}
\begin{minipage}[l]{0.5\textwidth}
 \centering \includegraphics[width=2in]{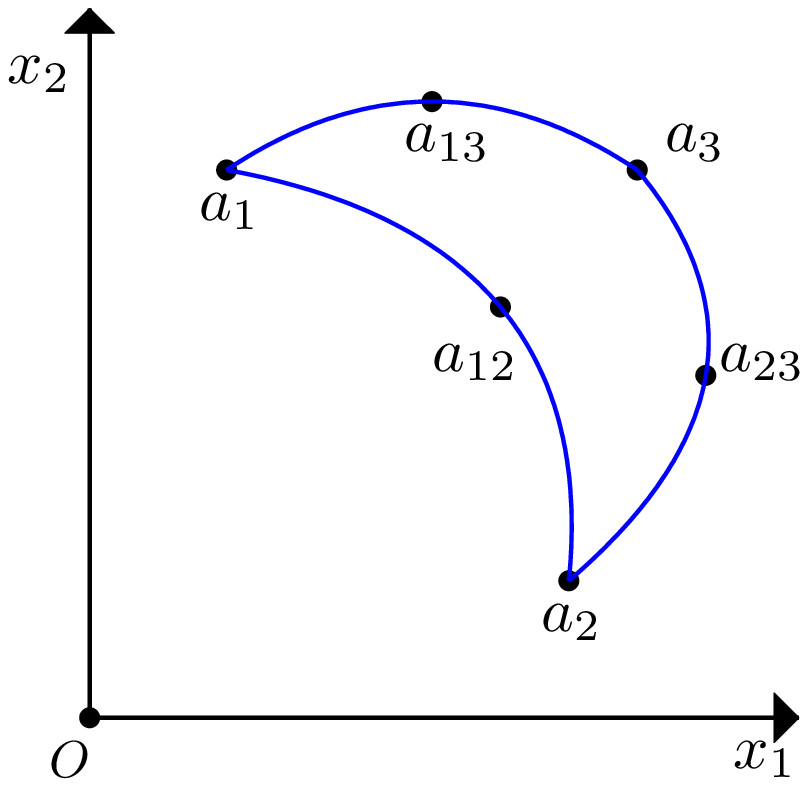}
 \caption{A curved triangular element T.}\label{dengcan curved element}
\end{minipage}
\end{figure}

Let $(\hat{T},\hat{P},\hat{\Sigma})$ be a quadratic Lagrange reference
element. Define $F_T:\hat{T}\rightarrow \mathbb{R}^2$
\begin{equation}\label{e:quadratic_mapping}
        \left\{
                \begin{aligned}
                   & F_T\in(P_2(\hat{T}))^2,\\
                   &x=F_T(\hat{x})=
                   \sum\limits_{i=1}^{3}a_i\hat{\mu}_i(\hat{x})+
        \sum\limits_{1\leq i<j\leq 3}a_{ij} \hat{\mu}_{ij}(\hat{x}),
                \end{aligned}
         \right .
    \end{equation}
where $a_i,1\leq i\leq 3$, and $a_{ij},1\leq i<j\leq 3$ are given points
in $\mathbb{R}^2$, and
$$
\hat{\mu}_i(\hat{x})=\hat{\lambda}_i(\hat{x})(2\hat{\lambda}_i
(\hat{x})-1), \quad \hat{\mu}_{ij}(\hat{x})=
4\hat{\lambda}_i(\hat{x})\hat{\lambda}_j(\hat{x}),
$$
with $\hat{\lambda}_i(\hat{x}),1\leq i\leq 3$ being the barycentric
coordinates of $\hat{T}$.
If the map $F_T$ defined above is an injection, then $T=F_T(\hat{T})$ is a
curved triangular element as shown in Figure \ref{dengcan curved element}.
The standard quadratic iso-parametric finite element is defined as a finite
element triple $(T,P_T,\Sigma_T)$ with
\begin{equation}
 \label{iso_element}
    \left\{
        \begin{aligned}
  &   T=F_T(\hat{T}) \text{ being a curved triangle element},\\
  &   P_T=\{p:T\rightarrow\mathbb{R}^2~|~p=\hat{p}\circ F_T^{-1},\;
  \hat{p}\in\hat{P}\},\\
  &   \Sigma_T=\{p(a_i),1\le i\le 3; p(a_{ij}),1\le i<j\le 3\}.
   \end{aligned}
   \right .
\end{equation}

\subsection{Large expansion accommodating triangulations}
    Let $\hat{\mathcal{J}}$ be a straight edged triangulation on
    $\Omega_{\varrho} = B_1(0) \setminus B_{\varrho}(0)$, and let $1>\mu > \varrho$ be given constants. For a triangular element $\hat{K} \in \hat{\mathcal{J}}$ with vertices
$a_i, 1 \leq i \leq 3$, to accommodate the large expansionary deformation
around the defect and approximate the curved boundary better, choose $a_{ij}$ in the following way: if $a_i$, $a_j \in \{x: \mu <|x|<1\}$, then $a_{ij}$ is set as the midpoint of $a_i$ and $a_j$; otherwise denote
$(r(x),\theta(x))$ the polar coordinates of $x$, then set
\begin{equation}
  \label{eq:make_curved_mesh}
  a_{ij}=(r_{ij}\cos\theta_{ij},r_{ij}\sin\theta_{ij}),
\end{equation}
where
\begin{equation*}
 r_{ij}=\frac{r(a_i)+r(a_j)}{2},\;\;\;
 \theta_{ij}=\frac{\theta(a_i)+\theta(a_j)}{2}.
 \end{equation*}
With the six points $a_i$ and $a_{ij}$, a triangular element
$K$ is defined by the mapping \eqref{e:quadratic_mapping}. With this kind
of curved elements in a neighborhood of the defects and on the outer boundary,
while using general straight triangles elsewhere, the
mesh can better accommodate the locally large expansionary deformations.

\begin{figure}[ht!]
     \begin{minipage}[l]{0.5\textwidth}
 \centering \includegraphics[width=2.4in]{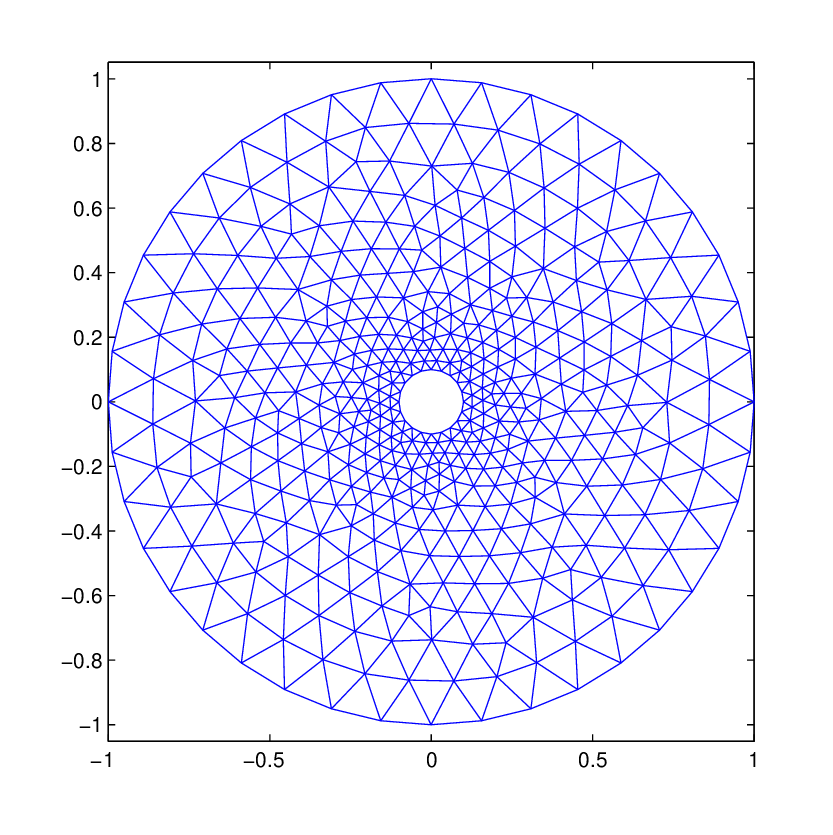}
\vspace*{-8mm}
 \caption{An EasyMesh $\hat{\mathcal{J}}'$.} \label{easy_mesh}
 \end{minipage}
\begin{minipage}[l]{0.5\textwidth}
 \centering \includegraphics[width=2.4in]{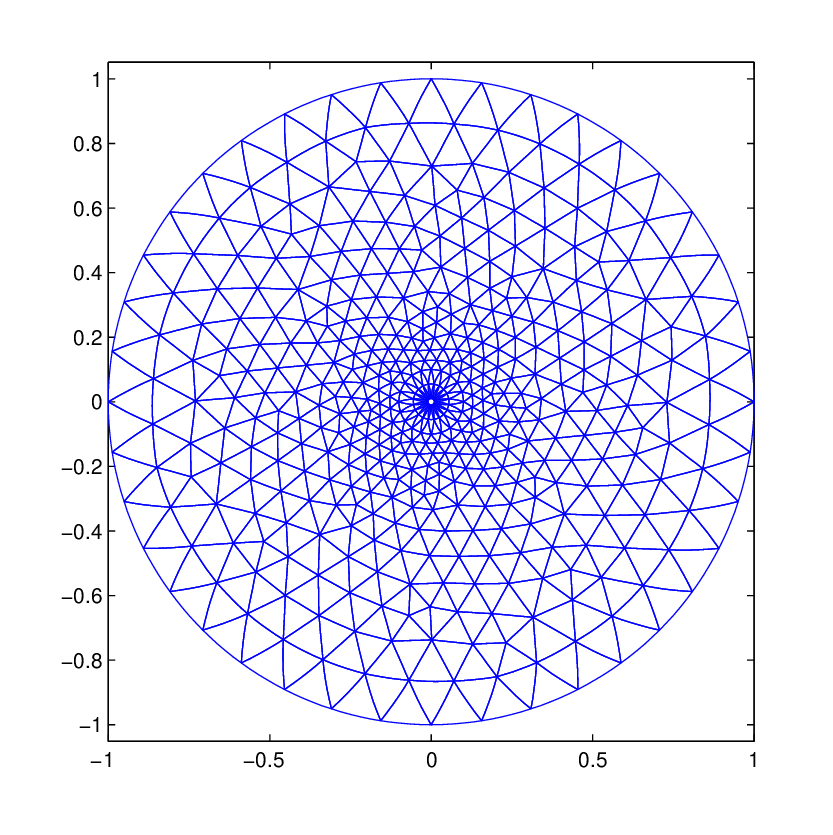}
\vspace*{-8mm}
 \caption{$\mathcal{J}$, a mesh adapted to cavity.} \label{final_mesh}
\end{minipage}
\end{figure}

When the radius of the defect $\varrho$ is very small, the mesh produced by
the EasyMesh (a software producing 2-d triangular mesh \cite{Lian and Li iso}) can be irregular near
the defect. To produce a mesh which can
better accommodate a
void growth, as suggested by \cite{Lian and Li iso}, a mesh
$\hat{\mathcal{J}}'$ can be introduced on $\Omega_{\hat{\varrho}}$
with $\hat{\varrho}\gg \varrho$ by the EasyMesh, which is then transformed as $\mathcal{J}'$
under the iso-parametric deformation as above and coupled with
several layers of circumferentially uniform mesh $\mathcal{J}''$ given on the domain
$\{x : \varrho \leq |x|  \leq \hat{\varrho} \}$, where on each layer the mesh is
similar to that shown in Figure~\ref{tuli} (also see therein for a more specific description),
to produce a final mesh $\mathcal{J} = \mathcal{J'} \cup \mathcal{J''}$. As an example,
an EasyMesh produced mesh $\hat{\mathcal{J}}'$ with $\hat{\varrho}=0.1$ is shown in
Figure~\ref{easy_mesh}, and the final mesh $\mathcal{J}$
with $\varrho=0.01$, $\mu=0.15$, and a two-layer
circumferentially uniform mesh $\mathcal{J}''$ is shown in Figure~\ref{final_mesh}.

\section{On the orientation-preservation conditions}

\begin{figure}[h!]

 \centering

 \begin{minipage}[c]{1\textwidth}

     \centering \includegraphics[width=7cm, height=7cm]{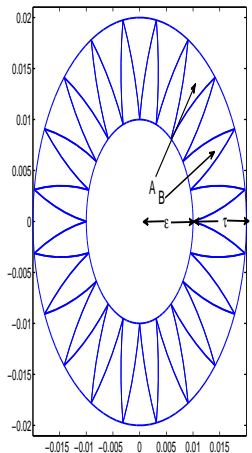}

 \vspace*{-2mm}

 \caption{A typical layer of circumferentially uniform curved triangulation.} \label{tuli}

 \end{minipage}

\end{figure}

We are concerned with orientation-preservation of large expansionary finite
element deformations around a
small prescribed void. A typical curved triangulation around a prescribed
circular ring with inner radius $\epsilon = 0.01$ and thickness $\tau =0.01$
is shown in Figure~\ref{tuli}, in which we see that the curved triangular
elements can be classified into two basic types, namely types A and B. More precisely,
let $N$ be the number of evenly spaced nodes on both circles, then, each of the
types A and B elements $A_i$, $B_i$, $i=0,1,2,\cdots,N-1$, are defined by three nodes, denoted as
$a^A_{i,j}=(r^A_{i,j},\theta^A_{i,j})$ and $a^B_{i,j}=(r^B_{i,j},\theta^B_{i,j})$,
$j=1,2,3$, in polar coordinates, as follows:
\begin{eqnarray*}
r^A_{i,1}=r^B_{i,2}=r^B_{i,3}=\epsilon, & & r^A_{i,2}=r^A_{i,3}=r^B_{i,1}=\epsilon + \tau, \\
\theta^A_{i,1}=\theta^B_{i,3}=\theta^B_{i-1,2}=\frac{2i}{N}\pi, & &
\theta^A_{i,2}=\theta^A_{i-1,3}=\theta^B_{i-1,1}=\frac{2i-1}{N}\pi,
\end{eqnarray*}
where the first subscript in $r$ and $\theta$ is understood in the sense of \!\!$\mod\!(N)$.

To have a better picture in our mind for the problem given below, we first introduce
some notations. Let $\epsilon$ and $\tau$ represent respectively the inner radius
and the thickness of the circular annulus as shown in Figure~~\ref{tuli}.
Let $2N$ be the number of the elements in the circular annulus, and denote
$\kappa \triangleq \epsilon/\tau$,
$\Omega_{(\epsilon,\tau)}=\{x \in \mathbb{R}^2: \epsilon \le |x| \le \epsilon+\tau\}$.
Throughout the paper, we use the notation
$\Phi \preceq \Psi$ to mean that there exists a generic constant $C$ independent of
$\epsilon$, $\tau$ such that $|\Phi| \le C \Psi$. And $\Phi\sim \Psi$ means that
$\Psi \preceq \Phi \preceq \Psi$.

Let $u$ be the cavitation solution, and let $\mathcal{J}=\mathcal{J}'\bigcup
\mathcal{J}''$ be a given mesh (see \S~2) with the layers in $\mathcal{J}''$ consisting
of well defined curved triangular elements of types A and B satisfying
$\det \nabla x>0$ (see Corollary~\ref{coro-mesh} for details).
To have $u$ well resolved by functions in the finite element function space
defined on $\mathcal{J}$, a necessary condition is that the finite element
interpolation function $\Pi u(x)$ is an admissible function. Since
$u$ is considered regular elsewhere other than in
$\mathcal{J}''$, where the material is subjected to a locally
large expansion dominant deformation, and since the key for a finite element
interpolation function to be admissible is $\det \nabla \Pi u(x)>0$ on
each of the curved triangular element, for simplicity and without loss of generality,
we will investigate in this section the conditions that ensure
$\det \nabla \Pi u(x)>0$ for radially symmetric expansionary
deformations of the form $u(x)=\frac{r(|x|)}{|x|}x$.
Since $\det \nabla \Pi u(x)\cdot\det \nabla x=
\det \frac{\partial \Pi u}{\partial \hat{x}}$, it suffices to ensure
$\det \frac{\partial \Pi u}{\partial \hat{x}}>0$ and
$\det \nabla x>0$ on the curved
triangular elements of types A and B in all layers (see \S~2.3).

The two lemmas below are the main ingredients for the orientation-preservation
conditions. To simplify the notations, for any positive function $s(\cdot)$, we denote
$s_0=s(\epsilon)$, $s_{1/2}=s(\epsilon+\tau/2)$, $s_1=s(\epsilon+\tau)$,
$\kappa^s_0=\frac{s(\epsilon)}{s(\epsilon+\tau)}$, and
$\kappa^s_{1/2}=\frac{s(\epsilon+\frac{\tau}{2})}{s(\epsilon+\tau)}$.

\vskip 2mm
We first give the orientation-preservation conditions for the type A elements.

\begin{lemma}\label{funA}
Let $v(x)=\frac{s(|x|)}{|x|}x$, where $s(t)$ is a
positive function satisfying $s'(t)> 0$, $s''(t)\ge 0$, $\forall$ $t \in (0, 1]$,
and
\begin{equation}\label{2r1/2>r1}
2s_{1/2} > s_1.
\end{equation}
Then, the Jacobian determinant
$\det \frac{\partial \Pi v(F_T(\hat{x}))}{\partial\hat{x}}$
of the iso-parametric finite element interpolation function
$\Pi v(F_T(\hat{x}))$
is positive on the elements of type A in $\Omega_{(\epsilon,\tau)}$ if and only if
 \begin{equation}\label{alp2}
-3s_0-s_1\cos{\frac{\pi}{N}}+4s_{1/2}\cos{\frac{\pi}{2N}}>0,
\end{equation}
 and
 \begin{equation}\label{ll}
-6s_1\cos^3{\frac{\pi}{2N}}+4s_{1/2}\cos^2{\frac{\pi}{2N}}+(s_0+9s_1)\cos{\frac{\pi}{2N}}-8s_{1/2}>0.
\end{equation}
\end{lemma}
\textbf{Proof.}
For the radially symmetric function
$v(x)=\frac{s(|x|)}{|x|} x$, the iso-parametric
finite element interpolation function can be written as (see \S~2)
\begin{equation}
\Pi v(x)=\sum\limits_{i=1}^{3}b_i\hat{\mu}_i(\hat{x})+
        \sum\limits_{1\leq i<j\leq 3}b_{ij}\hat{\mu}_{ij}(\hat{x}),
\label{eq deformation}
\end{equation}
where $\hat{x}=F_T^{-1}(x)$, and where, for
a representative of type A element,
$b_1=(s_0,0)$,
$b_2=(s_1\cos{\frac{\pi}{N}},-s_1\sin{\frac{\pi}{N}})$,
$b_3=(s_1\cos{\frac{\pi}{N}},s_1\sin{\frac{\pi}{N}})$,
$b_{12}=(s_{1/2}\cos{\frac{\pi}{2N}},-s_{1/2}\sin{\frac{\pi}{2N}})$,
$b_{13}=(s_{1/2}\cos{\frac{\pi}{2N}},s_{1/2}\sin{\frac{\pi}{2N}})$,
$b_{23}=(s_1,0)$.
On this element, we have
$$\Pi v(x)=(s_0+\alpha_2 y+2\alpha_1 y^2-4s_1\sin^2{\frac{\pi}{2N}}(\hat{x}_1^2+\hat{x}_2^2),
(2\gamma y-\beta)(\hat{x}_2-\hat{x}_1)),$$
where $y=\hat{x}_1+\hat{x}_2$,
\begin{eqnarray}
\alpha_1 &=&s_0+s_1-2s_{1/2}\cos{\frac{\pi}{2N}}, \label{alpha_1} \\
\alpha_2&=&-3s_0-s_1\cos{\frac{\pi}{N}}+4s_{1/2}
\cos{\frac{\pi}{2N}}, \label{alpha_2}\\
\beta&=&s_1 \sin{\frac{\pi}{N}}-4s_{1/2}
\sin{\frac{\pi}{2N}}, \label{beta} \\
\gamma&=&s_1 \sin{\frac{\pi}{N}}-2s_{1/2}
\sin{\frac{\pi}{2N}}.\label{gamma}
\end{eqnarray}
Hence
\begin{equation}\label{Jacobian Pi u}
\frac{\partial \Pi v}{\partial \hat{x}}=\left(
\begin{array}{cc}
\alpha-8s_1\hat{x}_1 \sin^2 {\frac{\pi}{2N}} & \alpha-8s_1\hat{x}_2 \sin^2
{\frac{\pi}{2N}} \\
\beta-4\gamma \hat{x}_1 & -\beta+4 \gamma \hat{x}_2 \\
\end{array}
\right),
\end{equation}
where $\alpha=4\alpha_1 y+\alpha_2$. It follows that
\begin{multline}\label{det Pi}
\det \frac{\partial \Pi v}{\partial \hat{x}}(\hat{x}_1,\hat{x}_2)=
H(y,z)\triangleq
16 \gamma \alpha_1 y^2-64s_1\gamma \sin^2{\frac{\pi}{2N}}z \\ +
\big(-8\beta \big(\alpha_1-s_1\sin^2{\frac{\pi}{2N}}\big)+
4 \gamma \alpha_2\big)y-2\beta \alpha_2,
\end{multline}
where $z=\hat{x}_1\hat{x}_2$.
Note that $\det \frac{\partial \Pi v}{\partial \hat{x}}>0$ on
$\hat{T}= \hat{T}_1 \cup \hat{T}_2$, where $\hat{T}_1 =
\{ (\hat{x}_1,\hat{x}_2) : 0\leq \hat{x}_2 \leq 1/2, \hat{x}_2 \leq \hat{x}_1
\leq 1-\hat{x}_2 \}$ and $\hat{T}_2= \{ (\hat{x}_1,\hat{x}_2) : 0\leq
\hat{x}_1 \leq 1/2, \hat{x}_1\leq \hat{x}_2 \leq 1-\hat{x}_1 \}$, is
equivalent to $H(y,z) > 0$ on the domain $\{0 \le y \le 1$, $0 \le z \le
\frac{y^2}{4}\}$.

Firstly, \eqref{alp2} follows from $H(0,0)=-2\beta \alpha_2>0$ and $\beta<0$, and \eqref{ll} is a direct result of $H(1,0)=4s_1\sin {\frac{\pi}{2N}}(-6s_1\cos^3{\frac{\pi}{2N}}+ 4s_{1/2}\cos^2{\frac{\pi}{2N}}+(s_0+9s_1)\cos{\frac{\pi}{2N}}-8s_{1/2})>0$.

Conversely, we infer from \eqref{alp2} and \eqref{ll} that $H(y,z)>0$ for $0\le y\le 1$, $0\le z \le \frac{y^2}{4}$ as following.
\begin{enumerate}
\item If $\cos{\frac{\pi}{2N}}> \kappa^s_{1/2}$, i.e., $\gamma>0$, then $\frac{\partial H}{\partial z}<0$. Thus,
it suffices to show $H(y,z) > 0$ on the curve $z=y^2/4(0 \le y \le 1)$.
On this curve, we have
\begin{eqnarray*}
H(y,z)&=&G(y) \triangleq 16 \gamma \alpha_3 y^2+(-8 \beta \alpha_3+4
\gamma \alpha_2)y-
2 \beta \alpha_2\\
&=&2(2\gamma y-\beta)(4\alpha_3 y+\alpha_2),
\end{eqnarray*}
where
$$
\alpha_3=s_1\cos^2{\frac{\pi}{2N}}-2s_{1/2}
\cos{\frac{\pi}{2N}}+s_0.$$
Consider the sign of $\alpha_3$. If $\alpha_3 > 0$, then by \eqref{alp2}, both
roots $y_1=\frac{\beta}{2\gamma}$ and
 $y_2=-\frac{\alpha_2}{4\alpha_3}$ of the equation $G(y)=0$ are negative.
So, it follows from $G''(y)=32\gamma\alpha_3>0$ that $G(y) > 0$ for $y>0$.

While if $\alpha_3 \leq 0$, we have $G(0)>0$, $G''(y)\leq 0$ and,
recalling that $s>0$, $s'> 0$ and $s''\ge 0$,
\begin{eqnarray*}
G(1)&=&2(2\gamma-\beta)\big(2s_1\cos^2{\frac{\pi}{2N}}-4s_{1/2}
\cos{\frac{\pi}{2N}}+s_0+s_1\big),\\ &\geq&2(2\gamma -\beta)
\big(2s_{1/2}\cos^2{\frac{\pi}{2N}}-4s_{1/2}\cos{\frac{\pi}{2N}}+2s_{1/2}\big)
\\ &=&4(2\gamma-\beta)s_{1/2}\big(\cos{\frac{\pi}{2N}}-1\big)^2>0,
\end{eqnarray*}
as a consequence, we infer that $G(y)>0$ on $[0,1]$.

\item If $\cos{\frac{\pi}{2N}}= \kappa^s_{1/2}$, i.e., $\gamma =0$, then,
$H(y,z)=-2\beta(4\alpha_3 y+\alpha_2)$ with $\beta <0$. A similar
argument as in (i) yields $H(0,z)=-2 \beta \alpha_2 > 0$, $H(1,z)=-2\beta(\alpha_2+4 \alpha_3) > 0$, hence
we conclude $H(y,z)$ is positive on the domain $\{0 \le y \le 1, 0 \le z \le
\frac{y^2}{4}\}$.

\item If $\cos{\frac{\pi}{2N}}< \kappa^s_{1/2}$, i.e., $\gamma <0$, and thus
$\frac{\partial H(y,z)}{\partial z}>0$, it is sufficient to guarantee
$H(y,z) > 0$ on the set $\{z=0, 0 \le y \le 1\}$.
On the curve $z=0$,
$$
H(y,0)=G(y)\triangleq 16 \gamma \alpha_1 y^2+(-8\beta \alpha_3 +
4\gamma \alpha_2) y-2\beta \alpha_2.
$$
Since $\gamma <0$, and $s>0$, $s''\ge 0$ implies $\alpha_1>0$, it follows that
$G''(y)<0$. Thus, $G(y)>0$ for $y\in [0,1]$ equivalents to $G(0)>0$
and $G(1)>0$.

By \eqref{alp2} and $\beta<0$, $G(0)>0$. On the other hand, it follows from \eqref{ll} that
$G(1)=4s_1\sin{\frac{\pi}{2N}}(-6s_1\cos^3{\frac{\pi}{2N}}+
4s_{1/2}\cos^2{\frac{\pi}{2N}}+(s_0+9s_1)\cos{\frac{\pi}{2N}}-8s_{1/2})>0$.
\end{enumerate}
Thus we are led to the conclusion. \hfill $\square$

Similarly, the sufficient and necessary condition for the elements of type B is as  follows.
\begin{lemma}\label{funB}
Under the same assumptions of Lemma \ref{funA}, the Jacobian determinant
$\det \frac{\partial \Pi v(F_T(\hat{x}))}{\partial\hat{x}}$
of the iso-parametric finite element interpolation function
$\Pi v(F_T(\hat{x}))$ is positive on the curved elements of type B in the circular domain $\Omega_{(\epsilon,\tau)}$ if and only if
\begin{equation}\label{lk}
2 s_0 \cos^2{\frac{\pi}{2N}}-4s_{1/2}
\cos{\frac{\pi}{2N}}+s_0+s_1< 0.
\end{equation}
\end{lemma}
\textbf{Proof.} Consider $\Pi v(x)$ defined by \eqref{eq deformation}
on a representative element of type B with
$b_1=(s_1,0)$, $b_2=(s_0\cos{\frac{\pi}{N}},s_0\sin{\frac{\pi}{N}})$,
$b_3=(s_0\cos{\frac{\pi}{N}},-s_0\sin{\frac{\pi}{N}})$,
$b_{12}=(s_{1/2}\cos{\frac{\pi}{2N}},s_{1/2}\sin{\frac{\pi}{2N}})$,
$b_{13}=(s_{1/2}\cos{\frac{\pi}{2N}},-s_{1/2} \sin{\frac{\pi}{2N}})$,
$b_{23}=(s_0,0)$. On this element, one has
$$\Pi v(x)=(s_1+\bar{\alpha}_2 y+2\bar{\alpha}_1 y^2-4s_0\sin^2{\frac{\pi}{2N}}(\hat{x}_1^2+\hat{x}_2^2), (2\bar{\gamma} y-\bar{\beta})(\hat{x}_1-\hat{x}_2)),$$
where
\begin{eqnarray*}
\bar{\alpha}_1&=&s_0+s_1-2s_{1/2}\cos{\frac{\pi}{2N}},\\
\bar{\alpha}_2&=&-3s_1-s_0\cos{\frac{\pi}{N}}+4s_{1/2}
\cos{\frac{\pi}{2N}}, \\
\bar{\beta}&=&s_0 \sin{\frac{\pi}{N}}-4s_{1/2} \sin{\frac{\pi}{2N}},
\\ \bar{\gamma}&=&s_0 \sin{\frac{\pi}{N}}-2s_{1/2} \sin{\frac{\pi}{2N}}.
\end{eqnarray*}
Hence
\begin{equation}
\frac{\partial \Pi v}{\partial \hat{x}}=\left(
\begin{array}{cc}
\bar{\alpha}-8s_0\hat{x}_1 \sin^2 {\frac{\pi}{2N}} & \bar{\alpha}-8s_0\hat{x}_2 \sin^2
{\frac{\pi}{2N}} \\
-\bar{\beta}+4\bar{\gamma} \hat{x}_1 & \bar{\beta}-4 \bar{\gamma} \hat{x}_2\\
\end{array}
\right),
\end{equation}
$$
\det \frac{\partial \Pi v}{\partial \hat{x}}=H(y,z)=
-16 \bar{\gamma} \bar{\alpha}_1 y^2 +64 s_0\bar{\gamma} \sin^2{\frac{\pi}{2N}}z+
(8\bar{\beta} \bar{\alpha}_3-4 \bar{\gamma} \bar{\alpha}_2)y+2
\bar{\beta} \bar{\alpha}_2,$$
where, recalling that $s>0$, $s'> 0$ and $s''\ge 0$, we get
\begin{eqnarray*}
    \bar{\alpha}_3 &=&\bar{\alpha}_1-s_0\sin^2{\frac{\pi}{2N}}\\
    &=&s_0\cos^2{\frac{\pi}{2N}}-2s_{1/2}\cos{\frac{\pi}{2N}}+s_1\\
&\ge & s_0\cos^2{\frac{\pi}{2N}}-
(s_0+s_1)\cos{\frac{\pi}{2N}}+s_1\\
&=&(\cos{\frac{\pi}{2N}}-1)(s_0\cos{\frac{\pi}{2N}}-s_1)>0.
\end{eqnarray*}
Since $\bar{\gamma}<0$, i.e., $\frac{\partial H(y,z)}{\partial z}<0$, thus, it suffices to guarantee $H(y,z) > 0$ on the curve $z=y^2/4$, for $0 \le y \le 1$.
On this curve, we have
\begin{eqnarray*}
H(y,y^2/4)&=& G(y)\triangleq -16\bar{\gamma} \bar{\alpha}_3y^2+(8\bar{\beta} \bar{\alpha}_3-4\bar{\alpha}_2
\bar{\gamma})y+2 \bar{\alpha}_2 \bar{\beta}\\
&=& -2(2 \bar{\gamma} y- \bar{\beta})(4 \bar{\alpha}_3 y+ \bar{\alpha}_2).
\end{eqnarray*}
Let $y_1, y_2$ be the two roots of $G(y)=0$. Since
$y_1=\frac{\bar{\beta}}{2\bar{\gamma}}>1$ and
$G''(y)=-32\bar{\gamma}\bar{\alpha}_3>0$, we see that $G(y)>0$ on $[0,1]$
equivalents to $y_2=\frac{-\bar{\alpha}_2}{4 \bar{\alpha}_3} > 1$, or
$$
4\bar{\alpha}_3+\bar{\alpha}_2=2 s_0 \cos^2{\frac{\pi}{2N}}-4s_{1/2}
\cos{\frac{\pi}{2N}}+s_0+s_1< 0.
$$
Hence, the proof is completed. \hfill $\square$

\begin{remark}
As is shown in the proof of Lemmas \ref{funA} and \ref{funB}, $\det \frac{\partial \Pi v}{\partial \hat{x}}>0$ is satisfied on the elements if and only if $\det \frac{\partial \Pi v}{\partial \hat{x}}>0$
on the three vertices of the type A elements, while on the midpoint of the inner circle edge
of the type B elements.
\end{remark}

\begin{theorem}\label{tt}
Under the assumptions of Lemma \ref{funA}, the Jacobian determinant
$\det \frac{\partial \Pi v(F_T(\hat{x}))}{\partial\hat{x}}$ is positive on the curved
elements in the circular domain $\Omega_{(\epsilon,\tau)}$ if and only if
\begin{equation}\label{asp2}
4s_{1/2} > 3s_0+s_1,
\end{equation}
and $\cos {\frac{\pi}{2N}}> \max \{l_1,l_2\}$, where $l_1$ is the smaller root of the equation
\begin{equation}\label{l1}
2s_0z^2-4s_{1/2}z+s_0+s_1=0,
\end{equation}
$l_2$ is the second root of the equation
\begin{equation}\label{l2}
-6s_1z^3+4s_{1/2}z^2+(s_0+9s_1)z-8s_{1/2}=0.
\end{equation}
\end{theorem}
\textbf{Proof.} Firstly, the inequality \eqref{asp2} implies that $l_1<1$ as well as the
bigger root of equation \eqref{l1} is greater than 1. On the other hand, if \eqref{asp2}
does not hold, then one has $l_1 \ge 1$, consequently there is no $N$ such that \eqref{lk}
is satisfied. Next, let $l_3$ be the smaller
root of the equation (see \eqref{alp2})
$$
2s_1z^2-4s_{1/2}z+3s_0-s_1=0,
$$
Then the inequalities \eqref{alp2}, \eqref{ll}, \eqref{lk} are satisfied if and only if
$\cos{\frac{\pi}{2N}}>\max\{l_1,l_2,l_3\}$. Substituting $l_3$ into
the left hand side of \eqref{l1}, we have
$$
2s_0l_3^2-4s_{1/2}l_3+s_0+s_1=2(s_0-s_1)(l_3^2-1)>0,
$$
which together with \eqref{asp2} implies $l_3<l_1$. Thus, the conclusion of the
theorem follows from Lemmas \ref{funA} and \ref{funB}. \hfill $\square$

The theorem allows us to work out an explicit condition for a mesh, defined on a
ring region with curved triangular elements of types A and B, to be
well defined in the sense that $\det \frac{\partial x}{\partial \hat{x}}>0$.

\begin{corollary}\label{coro-mesh}
Let $\epsilon>0$, $\tau>0$, and $\kappa \triangleq \epsilon/\tau$. Let a layer
of evenly spaced $N$ couples of curved triangular mesh elements of types A and B
be introduced on the circular domain $\Omega_{(\epsilon,\tau)}$ by
\eqref{e:quadratic_mapping} and \eqref{eq:make_curved_mesh} (see Figure~\ref{tuli}).
Then, there exists an integer $\hat{N}(\kappa)$,
such that, the Jacobian determinant of the mesh map is positive, i.e.,
$\det \frac{\partial x}{\partial \hat{x}} >0$, if and only if
$N \geq \hat{N}(\kappa)$. Moreover,
$\hat{N}=\hat{N}(\kappa)\sim 1+\kappa^{1/4}$.
\end{corollary}
\textbf{Proof.} Taking $v(x)=x$, or equivalently $s(t)=\text{id}(t)=t$, in
Lemma~\ref{funA}, then it is easily verified that  \eqref{2r1/2>r1} and \eqref{asp2}
are satisfied. Consequently, by Theorem \ref{tt}, we conclude that
$\det \frac{\partial x}{\partial \hat{x}}>0$ on the elements if and only if
$\cos{\frac{\pi}{2N}}>\hat{l}(\kappa)=\max \{\hat{l}_1,\hat{l}_2\}$,
where $\hat{l}_1$ is the smaller root of the equation
\begin{equation} \label{hl1}
\kappa z^2-(1+2\kappa)z+\kappa+1/2=0,
\end{equation}
$\hat{l}_2$ is the second root of the equation
\begin{equation} \label{hl2}
-6(1+\kappa)z^3 + 4(1/2 +\kappa)z^2 +(9+10\kappa)z -
8(1/2+\kappa)=0.
\end{equation}
Thus the conclusion follows by setting
$\hat{N} \triangleq \big[\frac{\pi}{2\arccos{\hat{l}(\kappa)}}\big]+1$.
What remains for us to show
now is $\hat{N}=\hat{N}(\kappa)\sim 1+\kappa^{1/4}$.

If $\kappa =\epsilon/\tau$ is bounded above by a constant $C \geq 1$, note that $\hat{l}
<\frac{\kappa+1/2}{\kappa+1}\leq \frac{1+2C}{2+2C}$, it follows that $\hat{N}\leq [\frac{\pi}{2\arccos{\frac{1+2C}{2+2C}}}] +1$.

Next we consider the case when $\kappa >C$. Notice that
$\cos{\frac{\pi}{2N}}> \hat{l}(\kappa)$ is equivalent to
$\sin^2{\frac{\pi}{4N}}< \frac{1-\hat{l}(\kappa)}{2}=
\min\{\frac{1-\hat{l}_1(\kappa)}{2},\frac{1-\hat{l}_2(\kappa)}{2}\}$, and for
$\kappa \gg 1$, we have
\begin{equation}\label{1-hl1}
\frac{1-\hat{l}_1(\kappa)}{2}=\frac{1}{2(1+\sqrt{1+2\kappa})}=
\frac{\kappa^{-1/2}}{2\sqrt{2}}+ O(\kappa^{-1}).
\end{equation}
On the other hand, since $\hat{l}_2(\kappa)$ is the second root of the equation \eqref{hl2},
then $\frac{1-\hat{l}_2(\kappa)}{2}$ is the second root of the equation
$$
-6(1+\kappa)(1-2z)^3+4(\kappa+1/2)(1-2z)^2+(9+10\kappa)(1-2z)-8(1/2+\kappa)=0.
$$
Denote $t=\frac{1}{1+\kappa}$, then, the equation can be rewritten as
$$
z^3-(7/6+t/6)z^2+\frac{5t}{24}z+\frac{t}{48}=0.
$$
By the root formula of a cubic equation(see \cite{cubic root}), its second root is given by
$$
\frac{1-\hat{l}_2(\kappa)}{2}=-\iota\cos{(\psi+\pi/3)}+\frac{7+t}{18},
$$
where
$\iota=2(-\frac{w}{3})^{1/2}$, $\cos{3\psi}=-\frac{q}{2}(-\frac{w}{3})^{-3/2}$,
$w=-\frac{49}{108}+\frac{17t}{216}-\frac{t^2}{108}$,
$q=-\frac{343}{2916}+\frac{25t}{486}+\frac{17t^2}{3888}-\frac{t^3}{2916}$.
Hence, by the Taylor expansion, one has
\begin{equation}\label{1-hl2}
    \frac{1-\hat{l}_2(\kappa)}{2}=\frac{t^{1/2}}{2\sqrt {14}}+O(t)=
\frac{\kappa^{-1/2}}{2\sqrt{14}}+O(\kappa^{-1}).
\end{equation}
Note that $\hat{N}(\kappa)=\big[\frac{\pi}{4\arcsin
\sqrt{\frac{{1-\hat{l}(\kappa)}}{2}}}\big]+1$, the proof is
completed by \eqref{1-hl1} and \eqref{1-hl2}. \hfill $\square$

\begin{lemma}\label{epsilon>tau^2}
Let $\varrho>0$ be such that the smooth cavity solution $r(\cdot)$ is well defined
on $[\varrho,1]$ and satisfies $r(R) \ge r_c>0$, \eqref{mR<r'<MR} and
$|r^{(3)}(R)|\le Q$. Then, there exists a constant $C>0$,
such that \eqref{asp2} holds, if $\epsilon \geq \max\{\varrho, C \tau^2\}$
and $\epsilon+\tau \leq 1$.
Furthermore, let
$\tilde{\tau}_0 = \min \big\{1-\epsilon,\sqrt{\frac{2r_c}{\max_{\xi \in
[\varrho,1]}r''(\xi) }}\big\}$, then \eqref{2r1/2>r1} holds
for all $\tau \in (0,\tilde{\tau}_0]$.
\end{lemma}
\textbf{Proof.} Taylor expanding $r(\epsilon+\frac{\tau}{2})$,
$r(\epsilon+\tau)$ at $\epsilon$, one gets
$$
4r_{1/2}-3r_0-r_1=r'(\epsilon)\tau+\frac{1}{12}
r^{(3)}(\xi_1)\tau^3-\frac{1}{6}r^{(3)}(\xi_2)\tau^3.
$$
Since $r'(\epsilon)\ge m\epsilon$ by \eqref{mR<r'<MR}, \eqref{asp2} follows by taking
$C= \frac{1}{2m} Q$.

On the other hand, by the Taylor expansion and $r''(x)>0$, we have
\begin{equation*}
2r_{1/2}-r_1>r(\epsilon)- r''(\eta)\frac{\tau^2}{2}, \;\;\text{for some} \;\; \eta \in (\epsilon,\epsilon+\tau).
\end{equation*}
This yields the inequality \eqref{2r1/2>r1}.
\hfill $\square$

\vskip 3mm
The orientation-preservation conditions on the mesh can now be given as follows,
where, to simplify the notations, we set
$\kappa_0\triangleq \kappa^r_0=\frac{r(\epsilon)}{r(\epsilon+\tau)}$,
$\kappa_{1/2}\triangleq \kappa^r_{1/2}=\frac{r(\epsilon+\frac{\tau}{2})}{r(\epsilon+\tau)}$.

\begin{theorem}\label{t1}
Let $u(x)=\frac{r(|x|)}{|x|}x$ be the cavity solution satisfying the conditions of
Lemma~\ref{epsilon>tau^2}, and $\Pi u(x)$ be the interpolation function
of $u(x)$ on a quadratic iso-parametric finite element function
space defined on a mesh consisting of only elements of types A and B.
Then, there exist constants $\tilde{\tau}$, $C>0$, and an integer
$\tilde{N}(\kappa,\kappa_0,\kappa_{1/2})$, such that
$\det \frac{\partial \Pi u}{\partial x}>0$ on each of the finite
elements, if the mesh satisfies the conditions that $\tau \leq \tilde{\tau}$,
$\epsilon \ge C \tau^2$ and $N\ge \tilde{N}(\kappa,\kappa_0,\kappa_{1/2})$.
Moreover, $\tilde{N}^{-1}(\kappa,\kappa_0,\kappa_{1/2})
\sim (\epsilon\tau)^{\frac{1}{4}}$.
\end{theorem}
\textbf{Proof.} Since $u(x)=\frac{r(|x|)}{|x|}x$ is the smooth minimizer of
\eqref{mini}, it follows from Lemma~\ref{r'r''>0} that Lemma~\ref{epsilon>tau^2} holds
for $r(R)$. Taking $v(x)=u(x)$, or equivalently $s(t)=r(t)$, in Lemma~\ref{funA},
let $\tilde{\tau}$, $C$ be given by Lemma~\ref{epsilon>tau^2}, then it follows from
Theorem~\ref{tt} and Corollary \ref{coro-mesh} that, on a mesh subject to the
constraints $\tau \leq \tilde{\tau}$ and $\epsilon \ge C \tau^2$, the interpolation
function $\Pi u$ is orientation preserving, i.e.,
$\det \frac{\partial \Pi u}{\partial x}>0$ if and only if
$N\ge \tilde{N} \triangleq \max \{\hat{N}, [\frac{\pi}{2\arccos l}]+1\}$,
where $l=\max \{l_1,l_2\}$ with $l_1$, $l_2$ being given in Theorem~\ref{tt} by
setting $s(t)=r(t)$. What remains to show is
$\tilde{N}^{-1}(\kappa,\kappa_0,\kappa_{1/2})\sim (\epsilon\tau)^{\frac{1}{4}}$.

Note that
\begin{eqnarray*}
1-l_1&=&\frac{\sqrt{r_{1/2}^2-r_0^2/2-r_0r_1/2}+r_0-r_{1/2}}{r_0}\\
&=&\frac{\sqrt{\frac{r'(\epsilon)}{2}\tau+\frac{r^{(3)}(\xi_1)}{24}\tau^3-
\frac{r^{(3)}(\xi_2)}{12}\tau^3+O(\epsilon^2\tau^2+\epsilon \tau^3+\tau^4)}
-O(\epsilon \tau+\tau^2)}{\sqrt{r_0}},
\end{eqnarray*}
where $\xi_1 \in (\epsilon,\epsilon+\tau/2)$, $\xi_2 \in (\epsilon, \epsilon+\tau)$.
Since $m \epsilon \le r'(\epsilon) \le M \epsilon$, $\epsilon \ge C \tau^2$,
then $1-l_1\preceq (\epsilon \tau)^{1/2}$. On the other hand, by taking
$C=\frac{Q}{2m}$ as in Lemma~\ref{epsilon>tau^2}, one has that
$\frac{r'(\epsilon)}{2}\tau+\frac{r^{(3)}(\xi_1)}{24}\tau^3-
\frac{r^{(3)}(\xi_2)}{12}\tau^3 \ge \frac{m\epsilon \tau}{2}-\frac{Q}{8}\tau^3\ge \frac{m}{4}\epsilon \tau$. Thus $1-l_1\sim (\epsilon \tau)^{1/2}$.
Denote $\tilde{l}_2=\frac{1}{2}(1-l_2)$, then $\tilde{l}_2$ is
the second root of the equation (see \eqref{l2})
\begin{equation}\label{siml2}
-6(1-2z)^3+4\kappa_{1/2}(1-2z)^2+(\kappa_0+9)(1-2z)-8\kappa_{1/2}=0.
\end{equation}
By the root formula of a cubic equation, $\tilde {l}_2$ is given by
$$
\tilde{l}_2=-\iota \cos{\big(\psi+\frac{\pi}{3}\big)}-\frac{\kappa_{1/2}}{9}+\frac{1}{2},
$$
where
$\iota=2(-\frac{w}{3})^{1/2}$, $\cos{3\psi}=-\frac{q}{2}(-\frac{w}{3})^{3/2}$,
$w=-\frac{3}{8}-\frac{\kappa_0}{24}-\frac{\kappa_{1/2}^2}{27}$,
$q=-\frac{\kappa_{1/2}}{8}+\frac{\kappa_0
\kappa_{1/2}}{216}+\frac{2\kappa_{1/2}^3}{729}$.
By the Taylor expansion, we obtain
$$
\iota=\frac{7}{9}\Big(1-\frac{17}{196}\frac{r'(\epsilon+\tau)}{r_1}\tau+
\frac{13r''(\epsilon+\tau)\tau^2}{392r_1}\Big)+O(\epsilon^2\tau^2+\tau^3),$$
$$\cos{3\psi}=1-\frac{243r'(\epsilon+\tau)\tau}{1372r_1}+O(\epsilon^2\tau^2+\tau^3).$$
This leads to
$$
\sin{3\psi}=\frac{9}{7}\sqrt{\frac{3r'(\epsilon+\tau)\tau}{14r_1}}\,\big(1+
O(\epsilon^2+\tau)\big),
$$
and
$$
\psi=\frac{3}{7}\sqrt{\frac{3r'(\epsilon+\tau)\tau}{14r_1}}\,\big(1+
O(\epsilon^2+\tau)\,\big).
$$
Consequently, we get
\begin{eqnarray*}
\tilde{l}_2&=&-\frac{\iota}{2} \cos{\psi}+\frac{\sqrt{3}}{2} \iota \sin \psi
-\frac{\kappa_{1/2}}{9}+\frac{1}{2}\\
&=&-\frac{\iota}{2}+\frac{\sqrt{3}}{2}\iota \psi-\frac{\kappa_{1/2}}{9}+\frac{1}{2}+O(\psi^2)\\
&=&\frac{1}{2}\sqrt{\frac{r'(\epsilon+\tau)\tau}{14r_1}}+O(\epsilon \tau+\tau^2),\\
&\sim &\sqrt{\epsilon \tau+\tau^2}.
\end{eqnarray*}
Hence, $1-l_2\sim ((\epsilon+\tau)\tau)^{1/2}$. The
conclusion of the theorem now follows by the definition of $\tilde{N}$,
Corollary \ref{coro-mesh} and $\arccos l=2\arcsin {\sqrt{\frac{1-l}{2}}}$.
\hfill $\square$

\begin{remark}\label{inter}
We would like to point out that the condition $\epsilon \ge C \tau^2$ on the mesh in
Theorem~\ref{t1} is not necessary. It is just a sufficient condition
to ensure \eqref{asp2}. For an incompressible cavity solution, \eqref{asp2} is in
fact unconditionally satisfied, and thus no restriction on the thickness $\tau$
is required. Of course this does not change the fact that the quadratic iso-parametric
element on its own is unstable for incompressible elasticity, however the result could 
be useful for a properly coupled mixed finite element method. In the proof above, we
could as well obtain $\tilde{N}^{-1}\preceq (\epsilon \tau+\tau^3)^{1/4}$ without
the condition $\tau \preceq \epsilon^{1/2}$.
\end{remark}

\begin{remark}
For the nonsymmetric cavitation deformation, under certain regularity assumptions on
the solution, we can apply similar methods as in \cite{SuLiRectan} to obtain a sufficient
condition for the interpolation function to be orientation preserving.
\end{remark}

\begin{figure}[h!]

 \centering

 \begin{minipage}[c]{1\textwidth}

     \centering \includegraphics[width=14.5cm, height=7.5cm]{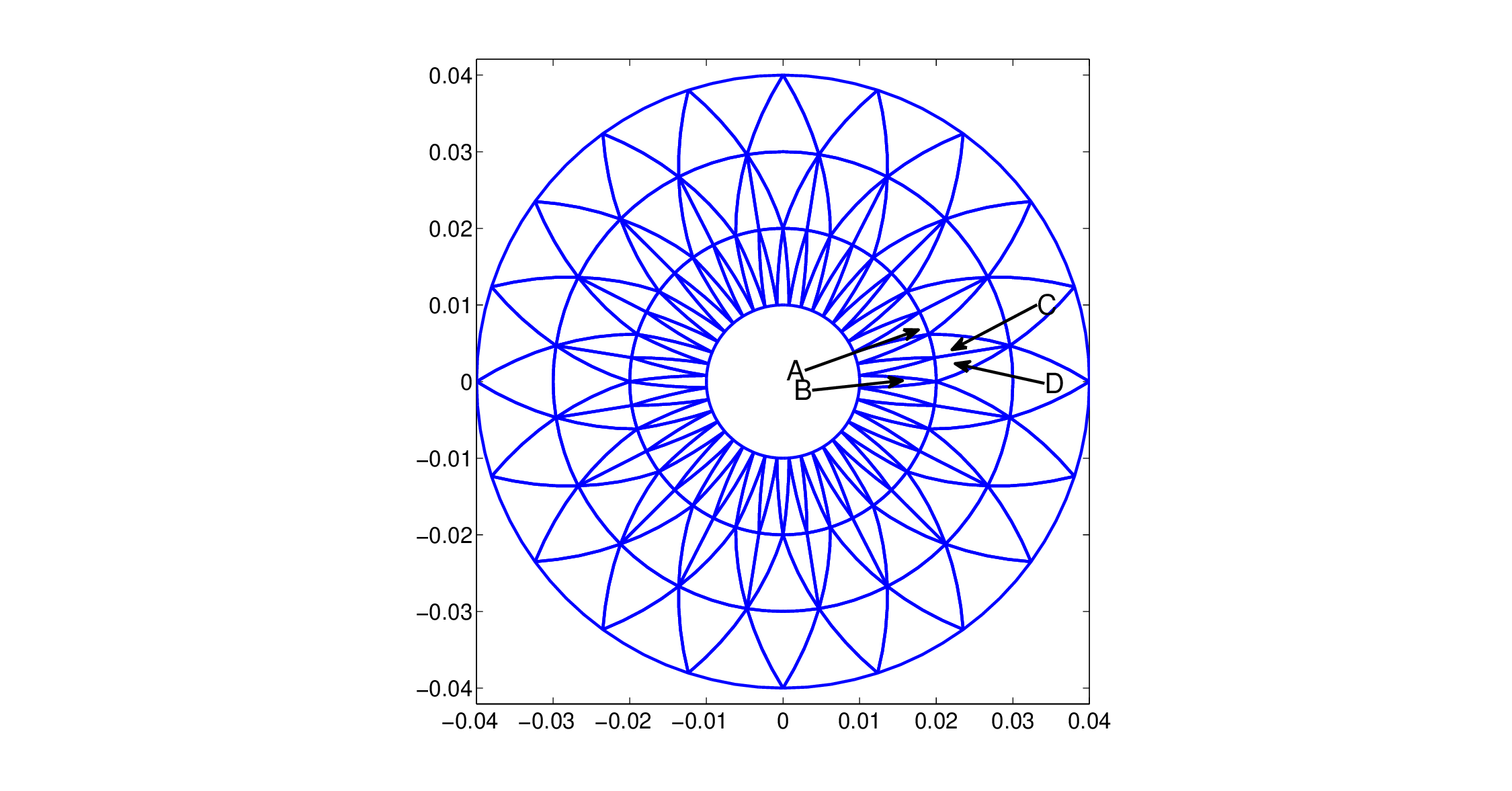}

 \vspace*{-2mm}

 \caption{Mesh coarsening is easily achieved.} \label{tuli1}

 \end{minipage}

\end{figure}

\begin{remark}
Compared to the orientation-preservation condition for the dual-parametric bi-quadratic
FEM in \cite{SuLiRectan}, where the corresponding sufficient and necessary condition for the
interpolation of the radially symmetric cavity solution is \eqref{asp2} only, while
the quadratic iso-parametric FEM imposes additional restrictions on the mesh
distribution in the angular direction, which can be more severe a condition. However,
when the radius $\varrho$ of the initial defect is very small,
to achieve the optimal interpolation error, similar restrictions on the mesh
distribution in the angular direction are also required for the dual-parametric
bi-quadratic FEM, particularly in the non-radially-symmetric case \cite{SuLi,SuLiRectan}.
Thus, to control the total degrees of freedom of the mesh, it is often necessary, for both
triangular and rectangular triangulations, to coarsening the mesh layers away from the cavity.
For our curved triangular partition, a conforming finite element mesh coarsening from
a circular ring layer to the next one outside can be easily achieved, by dividing each
type B element in the outside layer into two (types C and D as shown in Figure~\ref{tuli1})
with a straight line right in the middle along the radial direction without deteriorating
the orientation preservation and the approximation property, while it can be hardly done
for the curved rectangular one without introducing an intermediate layer.
Hence, for a conforming finite element cavity approximation, the quadratic
iso-parametric FEM can still be advantageous.
\end{remark}

As a comparison, we present below the orientation-preservation condition for the conforming
affine element.
\begin{theorem}\label{ta}
For the cavitation solution $u(x)=\frac{r(|x|)}{|x|}x$, the interpolation function in
the conforming affine finite element space is orientation preserving if and only if
$N \ge \tilde{N}_a$, where
$((\epsilon+\tau)\tau)^{-1/2}\preceq \tilde{N}_a \preceq (\epsilon \tau)^{-1/2}$.
\end{theorem}
\textbf{Proof.} For the radially symmetric deformation ~$v(x)=\frac{s(|x|)}{|x|}x$,
the interpolation function in the affine finite element space is given by
$\Pi v(x)=\sum\limits_{i=1}^3 b_i\hat{\lambda}_i(\hat{x})$. As in Figure \ref{reverse},
we can work on a typical triangle with $b_1=(0,s_0)$, $b_2=s_1(\sin {\frac{\pi}{N}},
\cos {\frac{\pi}{N}})$ and $b_3=s_1(-\sin {\frac{\pi}{N}},\cos {\frac{\pi}{N}})$. Thus
$$
\Pi v(x)=(s_1\sin {\frac{\pi}{N}}(\hat{x}_1-\hat{x}_2),s_0(1-\hat{x}_1-\hat{x}_2)
+s_1\cos {\frac{\pi}{N}}(\hat{x}_1+\hat{x}_2)),
$$
$\det \frac{\partial \Pi v}{\partial \hat{x}}=2s_1\sin{\frac{\pi}{N}}(s_1\cos
{\frac{\pi}{N}}-s_0)$. Hence the mesh is well defined, i.e.,
$\det \frac{\partial x}{\partial \hat{x}}>0$ if and only if
$\cos {\frac{\pi}{N}}>\frac{\epsilon}{\epsilon+\tau}$. It follows that for
$u(x)=\frac{r(|x|)}{|x|}x$, $\det \frac{\partial \Pi u}{\partial x}>0$ if and only
if $N>\tilde{N}_a \triangleq \max\{\frac{\pi}{\arccos{\frac{r(\epsilon)}{
r(\epsilon+\tau)}}},\frac{\pi}{\arccos{\frac{\epsilon}{\epsilon+\tau}}}\}$.
The conclusion is then established by
\eqref{mR<r'<MR}. \hfill $\square$

In \cite{Xu and Henao}, it is shown that a necessary condition for the conforming
piecewise affine finite element interpolation function of a cavity solution to have
finite energy in the layer $\Omega_{(\epsilon,\tau)}$ is $\tau \preceq
\epsilon^{p-1}$, where $p$ is the parameter in the energy density function
\eqref{energy function}. For $p=3/2$, this coincides with the condition
$\tau \preceq \epsilon^{1/2}$ used in Theorem~\ref{t1}. It is interesting to see
that, on a circular ring domain $\Omega_{(\epsilon,\tau)}$, by Theorems~\ref{t1}
and \ref{ta}, $\tilde{N}\sim \tilde{N}_a^{1/2}$ when $\tau \preceq \epsilon$, while
$\tilde{N}_a^{1/2} \preceq \tilde{N} \preceq \tilde{N}_a^{3/4}$ when
$\tau^2 \preceq \epsilon \preceq \tau$, i.e., the quadratic
iso-parametric finite element approximation needs significantly less elements.
For $2>p>3/2$, when the cavitation solution is harder to obtain numerically, the
restriction on the mesh for the conforming piecewise affine FEM is harsher,
which means a much larger number of total degrees of freedom is required. The fact, that the
number of elements needed on a layer with $\epsilon$ small so much
exceeds one's intuitive expectation, partially explains why no successful attempt
has ever been made at applying the affine FEM to the cavitation computation.

To illustrate the potential of our analysis in cavitation computation, we present below
some numerical results. The energy density in the numerical experiments is given by
\eqref{energy function} with $p=3/2$, $\omega=2/3$, and
$g(x)=2^{-1/4}(\frac12(x-1)^2+\frac{1}{x})$, the domain is
$\Omega_{0.01}\subseteq \mathbb{R}^2$ with a displacement boundary condition 
$u_0(x)=2x$ given on $\Gamma_0=\partial B_1(0)$ and a traction free boundary condition 
given on $\Gamma_1=\{x: |x|=0.01\}$.

\begin{figure}[h!]
 \centering
 \begin{minipage}[c]{0.68\textwidth}
     \centering \includegraphics[width=7cm, height=7cm]{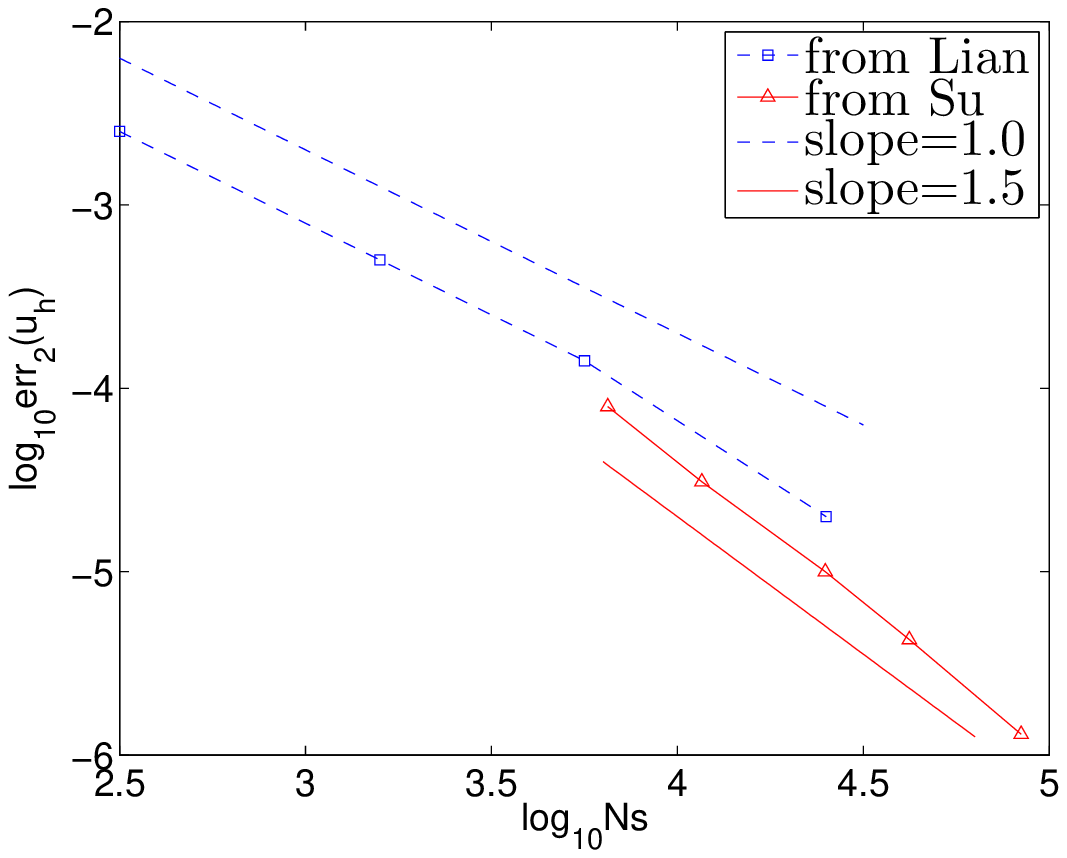}
 \vspace*{-2mm}
 \caption{$L^2$ error of numerical cavity solutions obtained on meshes based on
 experiences and a priori analysis.} \label{bijiao}
 \end{minipage}
\end{figure}

Figure \ref{bijiao} compares the $L^2$ error of the finite element cavity solutions $u_h$
against the total degrees of freedom $N_s$, 
where our result is obtained on the meshes produced according to our analytical results 
(near the cavity it is essentially governed by the orientation-preservation condition, 
see also \cite{SuLi}), while the meshes used in \cite{Lian and Li iso} were provided according to 
limited numerical experiences and thus, to guarantee the orientation preservation, the thickness
of the circular annulus were taken much thicker than necessary in general. 
It is clearly seen that our mesh is better in convergence rate
as well as actual accuracy.

\section{Conclusion remarks and discussions}

The orientation-preservation condition, i.e., the Jacobian determinant of the
deformation gradient $\det \nabla u>0$, is a natural physical
constraint in elasticity as well as in many other fields. It is well known
that the constraint can often cause serious difficulties in both theoretical analysis and
numerical computation, especially when the material is subject to large deformation as
in the case of cavitation. To overcome such difficulties can be crucial to successfully
solve the related problems.

In this paper, we analyzed the quadratic iso-parametric finite element
interpolation functions of the radially symmetric cavitation deformation on a class of
large radially symmetric expansion accommodating meshes, and obtained a set of sufficient
and necessary conditions on the orientation-preservation, which provide
a practical quantitative guide for the mesh distribution in the neighborhood of a cavity
in both radial and angular directions. Furthermore, the result shows that the
orientation-preserving cavitation approximation can be achieved by the quadratic
iso-parametric finite element method with a reasonable number of total degrees of
freedom, which is significantly smaller than the conforming piecewise affine finite
element method and is somehow comparable to the bi-quadratic dual-parametric
finite element method \cite{SuLiRectan}.
In fact, the orientation-preservation conditions together with the interpolation
error estimates, which will be established in a separate paper of ours \cite{SuLi},
will allow us to establish, for the quadratic iso-parametric FEM, a meshing strategy
leading to numerical cavitation solutions with optimal error bounds comparable to the
ones obtained in \cite{SuLiRectan} for a dual-parametric bi-quadratic FEM.

\end{document}